\documentclass[11pt,reqno]{amsart} 

\usepackage{textcomp}
\usepackage{amsmath}
\usepackage{amssymb} 
\usepackage[import,frame]{xy}
\usepackage{epsfig} 
\input ./warmread.sty 
\let\xyWARMprocess\xyWARMprocessMo 
\let\WARMprocessEPS\WARMprocessMoEPS 
\wrmquiet

\swapnumbers 
\theoremstyle{plain}
\newtheorem{theorem}{Theorem}[section]
\newtheorem{prop}[theorem]{Proposition}
\newtheorem{lemma}[theorem]{Lemma}
\newtheorem{cor}[theorem]{Corollary}

\theoremstyle{definition}
\newtheorem*{ack}{Acknowledgments}
\newtheorem{defn}[theorem]{Definition}

\newtheorem{example}[theorem]{Example}
\newtheorem{remark}[theorem]{Remark}
\newtheorem{notation}[theorem]{Notation} 
\newtheorem{question}[theorem]{Question}

  \def\m{\mathcal M_T}
  
\def\N{\mathbb N} \def\R{\mathbb R} \def\Z{\mathbb Z}

\def\phi{\varphi} \def\rho{\varrho}  \def\T{\mathbb
  T} \def\D{\mathcal D_F} \def\A{\mathcal A} 
\def\b{\mathcal B} \def\c{\mathcal C}\def\d{\mathcal D}
\def\lip{{\rm Lip}}

\begin{document}

\newcounter{algnum}
\newcounter{step}
\newtheorem{alg}{Algorithm}

\newenvironment{algorithm}{\begin{alg}\end{alg}}

\title[Flattening functions on flowers]
{Flattening functions on flowers}

\author{Edmund Harriss and Oliver Jenkinson}

\address{Edmund Harriss; School of Mathematical Sciences, Queen
  Mary, University of London, Mile End Road, London, E1 4NS, UK.
\newline {\tt edmund.harriss@mathematicians.org.uk} \newline {\tt
    www.mathematicians.org.uk/eoh}}

\address{Oliver Jenkinson; School of Mathematical Sciences, Queen
  Mary, University of London, Mile End Road, London, E1 4NS, UK.
  \newline {\tt omj@maths.qmul.ac.uk} \newline {\tt
    www.maths.qmul.ac.uk/$\sim$omj}}

\date\today
\maketitle

\begin{abstract}
Let $T$ be an orientation-preserving Lipschitz expanding map of the
circle $\T$. A pre-image selector is a map $\tau:\T\to\T$ with finitely
many discontinuities, each of which is a jump discontinuity, and such
that $\tau(x)\in T^{-1}(x)$ for all $x\in\T$.
The closure of the image of a pre-image selector is called a flower,
and a flower with $p$ connected components is called a $p$-flower.
We say that a real-valued Lipschitz function can be 
Lipschitz flattened on a 
flower whenever it is Lipschitz cohomologous to a constant on that flower.

The space of Lipschitz functions which can be flattened on a given
$p$-flower is shown to be of codimension $p$ in the space of all
Lipschitz functions, and the linear constraints determining this
subspace are derived explicitly.
If a Lipschitz function $f$ has a maximizing measure $S$ which is
Sturmian (i.e.~is carried by a $1$-flower), it is shown that $f$ can
be Lipschitz flattened on some $1$-flower carrying $S$.
\end{abstract}

\section{Introduction}
\label{introsection}

This article is concerned with certain special representatives of
dynamically defined cohomology classes, mainly motivated by the theory
of maximizing measures (see e.g.~\cite{bousch,bouschwalters,
  contreraslopesthieullen,lectures} for some background to this area).
Given a continuous self-map $T:X\to X$ of a compact metric space $X$,
let $\m$ denote the set of $T$-invariant Borel probability measures.
For a continuous function $f:X\to\R$, a measure $\mu\in\m$ is called
\emph{maximizing} if $\int f\, d\mu =\max_{m\in\m}\int f\, dm$.

One method of determining the maximizing measure(s) for a function $f$
consists of finding a continuous function $\varphi$ such that the set
$M(\tilde f):= \tilde f^{-1}(\max \tilde f)$ of maxima of $\tilde
f:=f+\varphi-\varphi\circ T$ carries\footnote{We say that a subset
  $G\subset X$ \emph{carries} a measure $\mu$ if the 
(topological) support of $\mu$,
 which we always denote by $\text{supp}(\mu)$, is contained in $G$.}  at least one
$T$-invariant measure.  The maximizing measures for $f$ (and for
$\tilde f$) are then precisely those invariant measures carried by
$M(\tilde f)$. It is known that such $\varphi$ always
exist if the map $T$
has some hyperbolicity and there is an appropriate control on the
modulus of continuity of the function $f$ (see e.g.~\cite{bousch,
  bouschwalters, conzeguivarch, contreraslopesthieullen,lectures}).

If $T$ is an expanding map, and $f$ is Lipschitz, the function
$\varphi$ may be chosen so that $T(M(\tilde f))=X$, i.e.~the set of
maxima of $\tilde f$ contains at least one pre-image of every point in
$X$ (see \cite[Thm.~1]{bouschwalters}). 
In this case it is easily seen that $M(\tilde f)$ carries at
least one invariant measure.  In certain 
situations it may be possible
to choose $\tilde f$ such that $M(\tilde f)$ 
contains \emph{precisely
  one} pre-image of every point in $X$; in this case we may think of
the function $f$
as determining a particular selection of pre-images of $T$.
The topology of the space $X$ may, however,
 preclude $M(\tilde f)$ from containing
exactly one pre-image of every point; for example this is the case if
$X=\T=\R/\Z$ is the circle, where no closed subset $F\subset \T$ is
mapped bijectively onto $\T$ by an expanding map $T$.  Nevertheless, 
it may
be the case that with \emph{finitely many exceptions}, every point in
$\T$ has a unique pre-image in $F=M(\tilde f)$.  For example if
$T:\T\to\T$ is defined by $T(x)=2x \pmod 1$ and $F$ is any closed
semi-circle, then apart from the common image of the endpoints of $F$,
every point in $\T$ has a unique pre-image in $F$.  More generally we
may choose $F$ to be the union of a finite number of disjoint
closed intervals with the property that the only points in $\T$
without a unique pre-image in $F$ are the images of the boundary
points in $F$. Such sets $F$ have been studied by Br\'emont
\cite{bremont1,bremont2}, who called them \emph{(finite)
  flowers} (some examples of flowers are depicted
in Figure~\ref{flowers} in Section \ref{sectiontwo}).

When $T(x)=2x \pmod 1$, Bousch \cite{bousch} has shown that if
$f:\T\to\R$ is any trigonometric polynomial of degree one then indeed
there exists $\tilde f=f +\varphi -\varphi\circ T$ such that $M(\tilde f)$ is a closed semi-circle (i.e.~$f$ satisfies the so-called 
\emph{Sturmian condition}, cf.~\cite{bousch}).  
Since closed
semi-circles have the notable property of 
carrying one, and only one,
$T$-invariant measure 
(so-called \emph{Sturmian} measures, see 
e.g.~\cite{bouschmairesse,bullettsentenac}), he deduced that every
degree-one trigonometric polynomial has a unique maximizing measure,
and that this measure is Sturmian.
The
analogous result has been obtained in \cite{adjr}
for a certain family of piecewise linear functions,
and experimental evidence suggests that Sturmian maximizing measures
appear rather often for sufficiently simple functions $f$.

If a Lipschitz function $f$ is known to have a Sturmian
maximizing measure\footnote{For example $f$ might be a degree-one
trigonometric polynomial, or a piecewise linear function as in
\cite{adjr}.} then 
a natural problem is
to determine precisely \emph{which}
of the Sturmian measures is $f$-maximizing.
We will establish (see Theorem \ref{flattentheorem})
the following 
useful necessary condition for 
a Sturmian measure $S$ to be $f$-maximizing:
there exists a closed semi-circle $F$ carrying $S$,
and a Lipschitz function $\varphi$, such that the restriction of
$\tilde f=f+\varphi-\varphi\circ T$ to $F$ is a constant function.
This condition was introduced by 
Bousch \cite{bousch}, who referred to
it as the \emph{pre-Sturmian} condition since it is clearly
implied by his Sturmian condition mentioned above.
The content of
Theorem \ref{flattentheorem} is that for Lipschitz $f$
whose maximizing measure is Sturmian yet which do \emph{not} satisfy
the Sturmian condition (such $f$ exist, 
cf.~Example~\ref{notmax}),
the pre-Sturmian condition is still satisfied
(on some semi-circle carrying the
maximizing measure). 

The practical utility of Theorem \ref{flattentheorem}
stems from the fact that, as noted by Bousch, 
solving the pre-Sturmian condition for $F$ amounts to
solving a \emph{real-valued} equation.
This real-valued equation,
which appears explicitly in \cite{bousch},
and as a particular case of our
Theorem \ref{equivalent}, can be solved numerically
(cf.~Remark \ref{zerofinding}).
It should be contrasted to the \emph{functional} equation
(namely, $\varphi(x)+(\max_{m\in\m}\int f\, dm)
=\max_{y\in T^{-1}(x)} (f+\varphi)(y)$, cf.~\cite{bousch})
which, a priori, must be solved for $\varphi$ in order to
check the Sturmian condition
(which in any case may not be
satisfied, cf.~Example \ref{notmax}).

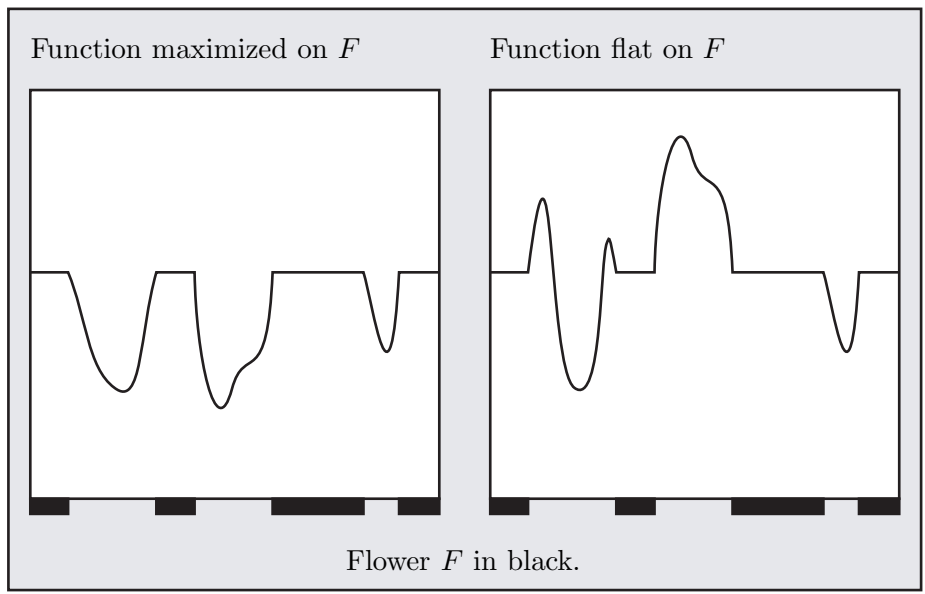
\begin{figure}[htpb]
$$\begin{xy}
\xyMarkedImport{}
        \xyMarkedPos{1}*!L\txt\labeltextstyle{Function maximized on $F$}
        \xyMarkedPos{2}*!L\txt\labeltextstyle{Function flat on $F$}
        \xyMarkedPos{3}*\txt\labeltextstyle{Flower $F$ in black.}
\end{xy}$$
\caption[Flattened functions]{Two functions flat on the flower $F$ 
(a $3$-flower for the doubling map $T_2(x)=2x \pmod 1$).  
The first function has the additional property that it attains
its global maximum precisely on $F$, so in particular is in normal form
(cf.~Definition \ref{normalformmax}); its maximizing measures 
are therefore precisely those invariant measures carried by $F$.}
\label{Flattened_functions}
\end{figure}

In the above discussion of the Sturmian and pre-Sturmian
conditions
the expanding map
 was assumed, for simplicity, to be
$T(x) =2x \pmod 1$.
In fact Theorem \ref{flattentheorem} 
is formulated in terms of arbitrary (orientation-preserving)
Lipschitz expanding maps of the circle.
In this general setting, closed semi-circles are replaced
by $1$-flowers
(i.e.~flowers with a single connected component).
Our interest in
the rest
of the article is in more general
\emph{$p$-flowers} (i.e.~flowers with $p$
connected components) for
orientation-preserving
Lipschitz expanding circle maps.
These are formally introduced in Section \ref{sectiontwo},
in terms of \emph{pre-image selectors} (roughly, a pre-image
selector is an inverse map for $T|_F$).
In Section \ref{sectionthree} we consider the generalization
of the pre-Sturmian condition
to arbitrary
$p$-flowers $F$:
for a Lipschitz function $f$,
if there exists a Lipschitz
$\varphi$ such that $(f+\varphi-\varphi\circ T)|_F$ is a constant
function,
we say that $f$ is \emph{Lipschitz flattened}
on $F$ (see 
Figure~\ref{Flattened_functions} in Section \ref{sectionthree}).
Theorem \ref{equivalent}
characterises this situation in terms of an
explicit vector-valued equation in $\R^p$. As with the real-valued
equation arising from the pre-Sturmian condition, this vector-valued
equation can in practice be solved  numerically in order to 
determine those $p$-flowers $F$ on which $f$ can be 
Lipschitz flattened.

Although nonlinear in $F$, the $p$ constraints detailed in
Theorem \ref{equivalent} depend linearly on the Lipschitz
function $f$. Theorem \ref{independentfunctionals} asserts
that these constraints are independent; in other words, the
space $\lip_F$ of Lipschitz functions which can be 
Lipschitz flattened
on a given $p$-flower $F$ is of codimension $p$ in the space
of all Lipschitz functions.

\begin{ack}
Both authors were partially supported by EPSRC grant
  GR/S50991/01.  The second author was partially supported by an EPSRC
 Advanced Research Fellowship.
\end{ack}

\section{Pre-image selectors and flowers}
\label{sectiontwo}

We consider the circle $\T$ as the space $\R/\Z$, i.e.~as the quotient
of the additive group of real numbers by the subgroup of integers.
The usual distance function on $\R$ induces a quotient distance
function on $\T$, which we shall denote by $d$.  The usual orientation
on $\R$ induces an orientation on $\T$.  If $a,b\in\T$ then $[a,b]$
denotes the positively oriented closed arc connecting $a$ to $b$, and
is called a \emph{closed interval}.  Open intervals $(a,b)$, and
half-open intervals $[a,b)$ and $(a,b]$ are defined analogously, in
accordance with the usual notational conventions.  We shall refer to
$a$ (respectively $b$) as the left (respectively right) endpoint of
any such interval.
For $a=c+\Z\in\T$, let $i:\T\to\R$ be the unique map
with image $[c,c+1)$ such that $i(u)\in u$ for all $u\in\T$.
Let $<_a$ denote the ordering on $\T$ induced by the usual 
ordering on $[c,c+1)$, i.e.~$u <_a v$ if and only if $i(u)<i(v)$.

A continuous map $T:\T\to\T$ is \emph{expanding} if there exists
$\varepsilon>0$ and $K>1$ such that for all $x,y\in\T$,
\begin{equation}
\label{expdef}
d(x,y)<\varepsilon\ \Rightarrow\ d(T(x),T(y))\ge Kd(x,y)\,.
\end{equation}

The degree of an expanding map is an integer of absolute value at
least 2 (see e.g.~\cite[p.~73]{katokhasselblatt}).  In this paper it
will be notationally convenient to only consider expanding maps which
are \emph{orientation-preserving}, i.e.~which have degree $k\ge 2$.
This, however, is not an essential restriction: results analogous to
those in this paper hold for orientation-reversing expanding maps, and
can be proved via slight modifications of the proofs given here.
Any expanding map of degree $k$ is topologically conjugate
to the map $T_k(x)= kx \pmod 1$
(see e.g.~\cite[p.~73, Thm.~2.4.6]{katokhasselblatt}); 
it will be useful to keep in mind
the maps $T_k$ as concrete examples for the notions and
results of this paper.

If the expanding map $T:\T\to\T$ has degree $k\ge 2$, and $a_0=T(a_0)$
is one of its $k-1$ fixed-points, then write
$T^{-1}(a_0)=\{a_0,\ldots,a_{k-1}\}$, where $a_0<\ldots<a_{k-1}$.
Define $X_{k-1}=[a_{k-1},a_0)$, and $X_i=[a_i,a_{i+1})$ for $0\le i\le
k-2$.  For $0\le i\le k-1$, the map
$$T_i:= (T|_{X_i})^{-1}$$
is called an \emph{inverse branch} of $T$.

We shall mainly be concerned with expanding maps $T$ which are in
addition \emph{Lipschitz} continuous, i.e.~there exists $C>0$ such
that
\begin{equation}
\label{lipschitzT}
d(Tx,Ty) \le C d(x,y)\quad\text{for all }x,y\in\T\,.
\end{equation}

Let $\text{Leb}$ denote normalised Lebesgue measure on $\T$, and let
$L^1$ denote the space of real-valued functions on $\T$ which are
integrable with respect to $\text{Leb}$.  For $g\in L^1$ we write
$\int g$ to denote the integral of $g$ with respect to $\text{Leb}$.
Let $L^\infty$ denote the space of functions $g:\T\to\R$ which are
essentially bounded with respect to $\text{Leb}$, and let
$\|g\|_{L^\infty}$ denote the essential supremum of $|g|$.
Let $L^\infty_0=\{g\in L^\infty : \int g=0\}$.
The following well-known result (see e.g.~\cite[Thm.~2.2.1]{ziemer})
will be used frequently.

\begin{lemma}
\label{rademachertheorem}
If $U:\T\to\T$ is piecewise Lipschitz, then its
derivative $U'$ exists Lebesgue almost everywhere, and defines an
element of $L^\infty$.  

If $U:\T\to\T$ is Lipschitz continuous
then $U'$ belongs to $L^\infty_0$.
Conversely, every
element of $L^\infty_0$ is the derivative of
a Lipschitz
continuous function on $\T$.
\end{lemma}

An orientation-preserving Lipschitz
map  $T:\T\to\T$ is easily seen to be expanding
if and only if there
exists $K>1$ such that $T'(x)\ge K$ for Lebesgue almost every
$x\in\T$.

\begin{figure}[htp]
$$\begin{xy}
        \xyMarkedImport{}
        \xyMarkedPos{1}*!L\txt\labeltextstyle{1-flower for $T_2$}
        \xyMarkedPos{2}*!L\txt\labeltextstyle{3-flower for $T_2$}
        \xyMarkedPos{3}*!L\txt\labeltextstyle{2-flower for $T_4$}
        \xyMarkedPos{4}*!L\txt\labeltextstyle{Pre-image selector}
        \xyMarkedPos{5}*!L\txt\labeltextstyle{Pre-image selector}
        \xyMarkedPos{6}*!L\txt\labeltextstyle{Pre-image selector}
        
        \xyMarkedPos{7}*\txt\labeltextstyle{\tiny{$0$}}
        \xyMarkedPos{8}*\txt\labeltextstyle{\tiny{$\alpha$}}
        \xyMarkedPos{49}*!U(.2)\txt\labeltextstyle{\tiny{$\alpha+$\textonehalf}}
        \xyMarkedPos{9}*\txt\labeltextstyle{\tiny{$0$}}
        \xyMarkedPos{10}*\txt\labeltextstyle{\tiny{$1$}}
        \xyMarkedPos{11}*!R\txt\labeltextstyle{\tiny{$1$}}
        \xyMarkedPos{12}*!R\txt\labeltextstyle{\tiny{$\alpha$}}
        \xyMarkedPos{13}*!R\txt\labeltextstyle{\tiny{$\alpha+$\textonehalf}}
        \xyMarkedPos{14}*\txt\labeltextstyle{\tiny{$2 \alpha$}}
        
        \xyMarkedPos{16}*\txt\labeltextstyle{\tiny{$\alpha$}}
        \xyMarkedPos{17}*\txt\labeltextstyle{\tiny{$\beta$}}
        \xyMarkedPos{18}*\txt\labeltextstyle{\tiny{$\gamma$}}
        \xyMarkedPos{50}*!U(.2)\txt\labeltextstyle{\tiny{$\alpha+$\textonehalf}}
        \xyMarkedPos{51}*!U(.2)!R(.3)\txt\labeltextstyle{\tiny{$\beta+$\textonehalf}}
        \xyMarkedPos{52}*!D(.7)!L(.5)\txt\labeltextstyle{\tiny{$\gamma+$\textonehalf}}
        \xyMarkedPos{22}*!R\txt\labeltextstyle{\tiny{$\gamma+$\textonehalf}}
        \xyMarkedPos{23}*!R\txt\labeltextstyle{\tiny{$\alpha+$\textonehalf}}
        \xyMarkedPos{24}*!R\txt\labeltextstyle{\tiny{$\beta+$\textonehalf}}
        \xyMarkedPos{25}*!R\txt\labeltextstyle{\tiny{$\gamma$}}
        \xyMarkedPos{26}*!R\txt\labeltextstyle{\tiny{$\alpha$}}
        \xyMarkedPos{27}*!R\txt\labeltextstyle{\tiny{$\beta$}}
        \xyMarkedPos{28}*\txt\labeltextstyle{\tiny{$2 \beta$}}
        \xyMarkedPos{29}*\txt\labeltextstyle{\tiny{$2 \alpha$}}
        \xyMarkedPos{30}*\txt\labeltextstyle{\tiny{$2 \gamma$}}
        
        \xyMarkedPos{34}*\txt\labeltextstyle{\tiny{$\alpha$}}
        \xyMarkedPos{36}*!U(.3)\txt\labeltextstyle{\tiny{$\alpha+$\textonequarter}}
        \xyMarkedPos{38}*\txt\labeltextstyle{\tiny{$\alpha+$\textonehalf}}
        \xyMarkedPos{32}*!L(.5)\txt\labeltextstyle{\tiny{$\alpha-$\textonequarter}}
        \xyMarkedPos{35}*\txt\labeltextstyle{\tiny{$\gamma$}}
        \xyMarkedPos{37}*\txt\labeltextstyle{\tiny{$\gamma+$\textonequarter}}
        \xyMarkedPos{39}*!D(.7)!L(.5)\txt\labeltextstyle{\tiny{$\gamma+$\textonehalf}}
        \xyMarkedPos{33}*\txt\labeltextstyle{\tiny{$\gamma-$\textonequarter}}
        \xyMarkedPos{43}*!R\txt\labeltextstyle{\tiny{$\gamma+$\textonehalf}}
        \xyMarkedPos{44}*!R\txt\labeltextstyle{\tiny{$\alpha+$\textonehalf}}
        \xyMarkedPos{45}*!R\txt\labeltextstyle{\tiny{$\alpha$}}
        \xyMarkedPos{46}*!R\txt\labeltextstyle{\tiny{$\gamma-$\textonequarter}}
        \xyMarkedPos{47}*\txt\labeltextstyle{\tiny{$4 \alpha - 1$}}
        \xyMarkedPos{48}*\txt\labeltextstyle{\tiny{$4 \gamma - 1$}}
\end{xy}$$
\caption[Flowers and pre-image selectors]{Flowers and their 
pre-image selectors (see
 Definition~\ref{defn:pre-image}).  
}
\label{flowers}
\end{figure}
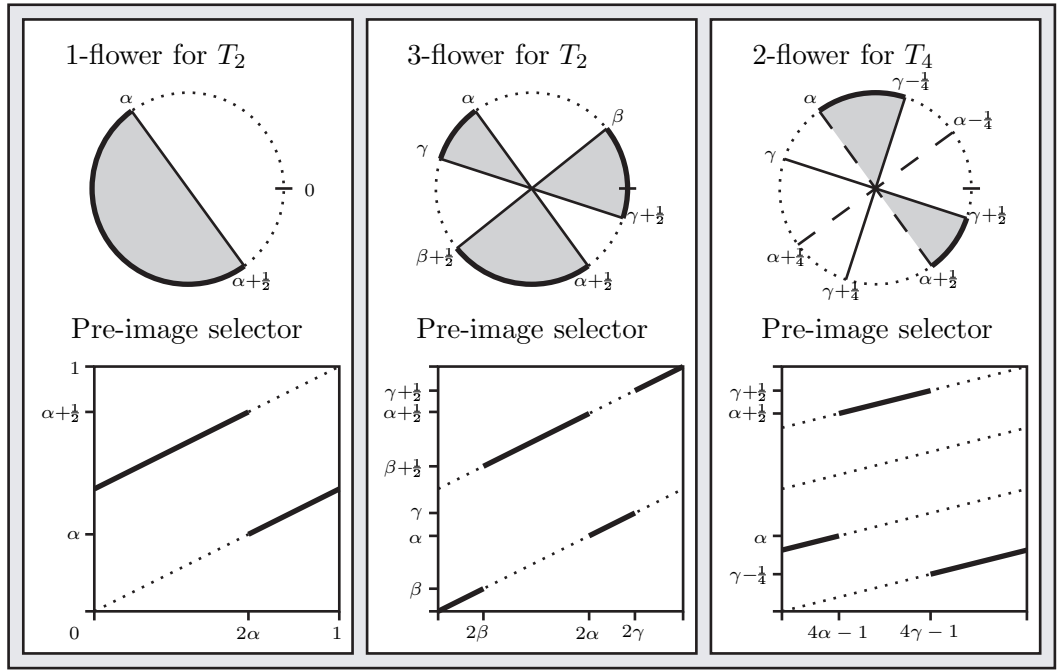

\begin{defn}
        \label{defn:pre-image}
  Let $T:\T\to\T$ be an expanding map.  A \emph{pre-image selector}
  for $T$ is a map $\tau:\T\to\T$ such that
\item[\, (i)]
$\tau(x)\in T^{-1}(x)$ for all $x\in\T$,
\item[\, (ii)] $\tau$ has finitely many discontinuities, each of which
  is a jump discontinuity (we say that $x$ is a jump discontinuity
  when $\lim_{z\nearrow x}\tau(z)$ and $\lim_{z\searrow x}\tau(z)$
  both exist but are distinct, and one of these values equals
  $\tau(x)$).
  
  If a pre-image selector $\tau$ has $p\ge1$ discontinuities then the
  closed set $F=\overline{\tau(\T)}$ is called the \emph{$p$-flower}
  (or simply the \emph{flower}) associated to $\tau$.  Each of the $p$
  connected components of $F$ is a closed interval with non-empty
  interior, and is called a \emph{petal} of $F$.

  If $F$ is a flower then define
  $$
  \D=T(\partial F)\,,
  $$
the set of discontinuities of any pre-image
  selector corresponding to $F$. Examples of flowers and pre-image selectors 
        are shown in Figure~\ref{flowers}.
\end{defn}

\WARMprocessEPS{flowers}{eps}{bb}

\begin{remark}\label{remark1}
\item[\, (a)]
  A pre-image selector $\tau:\T\to\tau(\T)$ is an inverse map for
  $T|_{\tau(\T)}$ (i.e.~$T\circ\tau$ is the identity on
  $\T$, and $\tau\circ T$ is the identity on $\tau(\T)$).  Since an
  expanding map on $\T$ is not a bijection, a pre-image selector is
  necessarily discontinuous.
\item[\, (b)]
In view of the discussion in Section \ref{introsection},
one might envisage relaxing the finiteness assumption in the 
definitions of pre-image selector and flower. The resulting objects
are more complicated than those considered here: 
for example the next simplest case would allow the
set of discontinuities of $\tau$ to be countably infinite. In this
case any accumulation point of the set of discontinuities would be a
discontinuity but not a jump discontinuity, and already this would
lead to certain complications in the proof of any analogue of
Theorem \ref{equivalent}.
\item[\, (c)]
We defined flowers above in terms of pre-image selectors.
Alternatively we could have defined flowers first, and then associated
a (non-unique) pre-image selector.
More precisely, if $F\subset\T$
can be written as
$F=\cup_{j=1}^p \overline{J_j}$, the closure of a
  finite union of open intervals $J_j$, such that $T(F)=\T$, and
  $T(J_i)\cap T(J_j)=\emptyset$ for $i\neq j$, then $F$ is the flower
  associated to some pre-image selector.  Indeed there are $2^p$ 
pre-image
  selectors $\tau$ whose associated flower is $F$. Each such $\tau$ is
  defined, on the \emph{interior} of $F$, 
by $\tau(x)=T^{-1}(x)\cap\left( \cup_{i=1}^p J_j\right)$ for
  $x\in \cup_{i=1}^p T(J_j)$; if $x\in T(\partial F)$ then
  $T^{-1}(x)\cap F=\{y,y'\}$ for some $y\neq y'$, and we may define
  $\tau(x)$ to equal either $y$ or $y'$.
\item[\, (d)]
Our terminology for flowers
differs slightly from that of
Br\'emont \cite{bremont1}.
What we call a \emph{$p$-flower} for a degree-$k$ map 
would, in \cite{bremont1}, be termed a $k$-flower with $p$ petals.

\item[\, (e)]
Every flower carries
 at least one $T$-invariant probability measure.
This is easily proved using the compactness of the flower,
and the fact that it 
contains a pre-image of each point on the
  circle.
\end{remark}

\begin{defn}
  Let $T:\T\to\T$ be an
orientation-preserving
expanding map, and $\tau$ a pre-image selector
  for $T$.  If $x$ is a discontinuity of $\tau$ then it is a jump
  discontinuity, so there exists $\varepsilon>0$ and
  $i,j\in\{0,\ldots,k-1\}$, $i\neq j$, such that
\begin{equation*}
\tau(z)=
\begin{cases}
  &T_i(z)\text{ for }x-\varepsilon<z<x\cr &T_j(z)\text{ for
  }x<z<x+\varepsilon\,,\cr
\end{cases}
\end{equation*}
and $\tau(x)$ equals either $T_i(x)$ or $T_j(x)$.  In this case we say
that $x$ is a \emph{discontinuity of type $(i,j)$}.  Note that both
the points $T_i(x)$ and $T_j(x)$ belong to the boundary $\partial F$
of the associated flower $F=\overline{\tau(\T)}$, and that no other
points in $T^{-1}(x)$ belong to $\partial F$.

For a discontinuity $x$ of $\tau$ we define $y=y(x)$ and $y'=y'(x)$ by
$$
y(x)= \lim_{z\searrow x}\tau(z)\quad,\quad y'(x)=\lim_{z\nearrow
  x}\tau(z)\ ,
$$
so in particular $y(x)$ is the left\footnote{If $T$ were
  orientation-\emph{reversing} then $y(x)$ would be a right endpoint,
  and $y'(x)$ a left endpoint.}  endpoint of a petal of $F$, and
$y'(x)$ is the right endpoint of a petal of $F$.\footnote{In general
  $y(x)$ and $y'(x)$ may, or may not be,
endpoints of the \emph{same} petal.}
\end{defn}

\begin{remark}
  A pre-image selector for an orientation-preserving Lipschitz
  expanding map is itself Lipschitz when restricted to any of its
 finitely many intervals of continuity.  
Therefore it is differentiable Lebesgue
  almost everywhere, and its derivative is in $L^\infty$,
by Lemma \ref{rademachertheorem}.  If the
  expanding map $T$ satisfies $T'(x)\ge K>1$ for Lebesgue almost every
  $x\in\T$, then the pre-image selector $\tau$ satisfies
\begin{equation}
\label{tauderivativebound}
0<\tau'(x)\le K^{-1}\quad\text{for Lebesgue almost every }x\in\T\,.
\end{equation}
In particular, the chain rule implies that
\begin{equation}
\label{iteratederivbound}
\|(\tau^n)'\|_{L^\infty}\le K^{-n}
\quad\text{for all }n\ge 1\,.
\end{equation}
\end{remark}

\begin{notation}
\label{orderingdefn}
  Let $\tau:\T\to\T$ be a pre-image selector, with corresponding
  flower $F$.  
Let $x_1$ denote that discontinuity
of $\tau$ which,
 with respect to the ordering $<_0$, is smaller
than all other discontinuities of $\tau$.\footnote{This choice of $x_1$
is for definiteness; in fact it is possible to fix 
$x_1$ to be an arbitrary
discontinuity of $\tau$.}
  Let $y_1$ be the unique pre-image of $x_1$ which is a left endpoint
  of some petal in $F$,
and for any $x\in\D$ define
\begin{equation*}
I_x=
\begin{cases}
  [y(x),y'(x)] & \text{if }y(x)<_{y_1} y'(x)\,,\cr [y'(x),y(x)] & \text{if
  }y'(x)<_{y_1} y(x)\,.
\end{cases}
\end{equation*}

Define
$$\A=\{x\in\D: I_x=[y(x),y'(x)]\}\quad,\quad \A'=\D\setminus \A\,.$$
\end{notation}

\begin{notation}
For a subset $G\subset \T$,
let $\chi(G)$ denote 
its characteristic
function.
\end{notation}

The following 
Lemma \ref{characteristicequality}
 will be a
very useful tool
in Section \ref{sectionthree}.
For maps $T$ of degree $k>2$, the combinatorics involved
in describing which pre-image of a point lies in a flower $F$ is
significantly more complicated than in the degree-2 case.
Lemma \ref{characteristicequality}, which is used
in the proof of
both Theorem \ref{equivalent}
and Theorem \ref{independentfunctionals},
 allows us to efficiently
sidestep these complications.

\begin{lemma}
\label{characteristicequality}
For any flower $F$,
$$\chi(F)=\sum_{x\in \A} \chi(I_x) - \sum_{x\in \A'} \chi(I_x)\,.$$
\end{lemma}
\begin{proof}
  Define $g:\T\to\R$ by
  $$g=\sum_{x\in \A} \chi(I_x) - \sum_{x\in \A'} \chi(I_x)\,.$$
  The
  function $g$ is upper semi-continuous and piecewise constant: its
  discontinuities are at the points $y(x)$ and $y'(x)$, for $x\in\D$.
  
  We claim that at every left endpoint $y$ of a petal of $F$, the
  function $g$ increases by $1$, in the sense that
  $g(y)=1+\lim_{z\nearrow y} g(z)$.  To see this, recall that every
  such left endpoint is of the form $y=y(x)$ for some $x\in\D$.  If
  $x\in\A$ then $\chi(I_x)=\chi([y(x),y'(x)])$ increases by $1$ at the
  point $y(x)$, while $g-\chi(I_x)$ is locally constant at $y(x)$, so
  $g=\chi(I_x)+(g-\chi(I_x))$ increases by $1$ at $y(x)$.  If
  $x\in\A'$ then $\chi(I_x)=\chi([y'(x),y(x)])$ decreases by $1$ at
  the point $y(x)$, while $g+\chi(I_x)$ is locally constant at $y(x)$,
  so $g=-\chi(I_x)+(g+\chi(I_x))$ increases by $1$ at $y(x)$.
  
  Similarly we can show that at every right endpoint $y'$ of a petal
  of $F$, the function $g$ decreases by $1$, in the sense that
  $g(y')=1+\lim_{z\searrow y} g(z)$.
  
  Now $y_1$ is the smallest point in the ordered set $(\T,<_{y_1})$, 
so $y_1\in
  I_{x_1}$, but $y_1\notin I_x$ for $x\in\D\setminus\{x_1\}$.
  Therefore $g(y_1)=1$.  With respect to $<_{y_1}$, the left and
  right endpoints of petals of $F$ alternate around the circle,
  beginning with a right endpoint after $y_1$.  Therefore $g$ takes
  the value $1$ between left and right endpoints of petals of $F$, and
  takes the value $0$ between right and left endpoints of petals of
  $F$.  That is, $g=\chi(F)$.
\end{proof}

\section{Lipschitz flattening}
\label{sectionthree}

Recall that 
$\m$ denotes the set of $T$-invariant Borel probability measures.
A continuous function $g:\T\to\R$ is called a \emph{weak coboundary}
if $\int g\, d\mu=0$ for every $\mu\in\m$, and a \emph{coboundary} if
$g=\varphi-\varphi\circ T$ for some continuous $\varphi:\T\to\R$.
Clearly every coboundary is a weak coboundary, but the converse is not
true (see e.g.~\cite{bouschjenkinson}).  However if $g$ is
\emph{Lipschitz}, then it is a coboundary if and only if it is a weak
coboundary, and if so then there is a \emph{Lipschitz} $\varphi$
(which is unique up to an additive constant) such that
$g=\varphi-\varphi\circ T$ (see \cite{livsic}).

\begin{defn}\label{flatteneddefn}
  Let $T:\T\to\T$ be an expanding map, and $F$ a corresponding flower.
A continuous function $f:\T\to\R$ is said to be \emph{flat} on $F$
if the restriction $f|_F$ is a constant function.
  
  We say that $f$ can be \emph{continuously flattened} on $F$ if there
  exists a weak coboundary $g:\T\to\R$ such that $f+g$ 
is flat on $F$.
  
  We say that $f$ can be \emph{Lipschitz flattened} on $F$ if there
  exists a Lipschitz coboundary $g:\T\to\R$ such that $f+g$
is flat on $F$ (or, equivalently, there exists a Lipschitz
  function $\varphi:\T\to\R$ such that $f+\varphi-\varphi\circ T$
  is flat on $F$).
\end{defn}

\WARMprocessEPS{Flattened_functions}{eps}{bb}

\begin{remark}
  If $f$ can be continuously flattened on $F$ then the constant
  function $(f+g)|_F$ is identically equal to $\int f\, d\mu$, where
  $\mu$ is any $T$-invariant probability measure carried by
  $F$.  If there are several such measures then this provides an
  obstruction
to being able to continuously flatten a
  given function on $F$: its integral must be the same with respect to
  each measure.  For example if $T(x)=2x \pmod 1$, and
  $F=[-\frac{1}{12},\frac{1}{12}]
  \cup[\frac{1}{4},\frac{5}{12}]\cup[\frac{7}{12},\frac{3}{4}]$, then
  $F$ contains both the fixed point $0$ and the period-$2$ orbit
  $\{1/3,2/3\}$, so a necessary condition for continuously flattening
  $f$ on $F$ is that $f(0)=\frac{1}{2}(f(1/3)+f(2/3))$. 
  \end{remark}

\begin{notation}
Let $\lip$ denote the vector space of all real-valued Lipschitz
functions on $\T$.
For a flower $F\subset\T$, define
$$
\lip_F=\{f\in\lip:f\text{ can be Lipschitz flattened on }F\}\,.
$$
\end{notation}

Using the various definitions and notation
introduced so far,
it is now possible to reiterate
the main results of this paper, in more detail than was possible
in Section \ref{introsection}.
The following is proved as
Theorems \ref{equivalent} and 
\ref{independentfunctionals}:

\smallskip

{\bf\noindent Theorem.} 
For any $p$-flower $F$, the set $\lip_F$ is a codimension-$p$
subspace of $\lip$.
Indeed there exist $p$ measures on $\T$, each absolutely continuous with
respect to Lebesgue measure, such that $\lip_F$ consists precisely
of those Lipschitz functions whose derivative has zero integral
with respect to each of these measures.

\smallskip

The Radon-Nikodym derivative of each of the $p$ measures 
in the above theorem 
is a certain infinite sum of characteristic functions of intervals
(explicit expressions
are given in Definition \ref{exdefnbelow} below,
see (\ref{exstar})).
These Radon-Nikodym derivatives can be rapidly approximated
(e.g.~by truncation of (\ref{exstar})), hence so can the kernels of the
corresponding measures (considered as functionals on $\lip$),
and therefore so can the
members of $\lip_F$.

As mentioned in Section \ref{introsection}, if a function $f$ can be
Lipschitz flattened on $F$, then the $T$-invariant measures carried by $F$
may, in certain circumstances, be good candidates for $f$-maximizing measures.
The relation between 
flattening and maximizing is particularly close when $F$ is a $1$-flower:
for any Lipschitz function $f:\T\to\R$ we have the following result,
which is proved later as Theorem \ref{flattentheorem}:

\smallskip

{\bf\noindent Theorem.} 
If the $f$-maximizing measure $S$ is carried by
some $1$-flower $F$, then $f$ can be Lipschitz flattened on some
$1$-flower $F'$ which carries $S$.

\smallskip

It should be noted that 
$F=F'$ for  non-atomic $S$, while
if $S$ is atomic
then $F'$ need not equal $F$.
The importance of the above theorem stems from the fact that invariant measures
carried by $1$-flowers (so-called \emph{Sturmian} measures,
cf.~Section 1 and Lemma \ref{sturmianfacts}) tend to arise
as $f$-maximizing measures for sufficiently simple functions $f$
(cf.~\cite{adjr,bousch,jenkinson2,major,announce}).
The theorem implies that for such $f$, the problem
of precisely identifying the
$f$-maximizing measure is reduced to determining those $1$-flowers
on which $f$ can be flattened, a problem which is computationally
accessible.

\begin{defn}\label{exdefnbelow}
  Let $\tau:\T\to\T$ be a pre-image selector, with corresponding flower
  $F$.  For each $x\in\D$, define $e_x:\T\to\R$ by
  \begin{equation}\label{exstar}
  e_x=\sum_{n=0}^\infty \chi(\tau^n I_x)\,.
  \end{equation}
  Note that $e_x\in L^1$, because the Lebesgue measure of $\tau^nI_x$
  decreases exponentially with $n$, by (\ref{iteratederivbound}).
\end{defn}

\begin{remark}
\label{escaperemark}
If $F$ is a \emph{1-flower} then a corresponding pre-image selector has
just one point of discontinuity $x$, and $F=[y,y']$, where
$y'=T_i(x)$, $y=T_{i+1}(x)$ for some $0\le i\le k-2$ (or
$y'=T_{k-1}(x)$, $y=T_0(x)$).  In this case the function $e_x$
will be denoted by $e_F$, and can be interpreted 
as the \emph{escape time function} for $F$,
$$e_x(t)=e_F(t)=\inf\{n\ge 0: T^n(t)\not\in F\}\,.$$
\end{remark}

If $f$ is itself Lipschitz then the following result, which
generalises \cite[Prop, p.~503]{bousch}, gives
  necessary and sufficient conditions for being able to
Lipschitz flatten $f$ on flowers.

\begin{theorem}
\label{equivalent}
Let $T:\T\to\T$ be an orientation-preserving Lipschitz 
expanding map.
Let $F\subset\T$ be a flower, and $\tau$ a corresponding 
pre-image selector.  If $f:\T\to\R$ is Lipschitz, then the following are
equivalent:
\item[\, (a)] $f$ can be Lipschitz flattened on $F$,
\item[\, (b)] For each discontinuity $x$ of $\tau$,
\begin{equation*}
\label{vanish1}
\int_{I_x} \sum_{n=0}^\infty (f\circ \tau^n)' =0\,.
\end{equation*}
\item[\, (c)] For each discontinuity $x$ of $\tau$,
\begin{equation*}
\int e_x f' = 0\,.
\end{equation*}
\end{theorem}
\begin{proof}
  (a) $\Rightarrow$ (b): If $f$ can be Lipschitz flattened on $F$ then
  there is a Lipschitz function $\varphi:\T\to\R$, and a constant
  $c\in\R$, such that $f+\varphi-\varphi\circ T=c$ on $F$.  In
  particular,
  $(f+\varphi-\varphi\circ T)(z)=c$
for all $z\in\tau(\T)$.
  Writing $z=\tau(x)$, and recalling that $\tau:\T\to\tau(\T)$ is a
  bijection such that $T\circ\tau$ is the identity on $\T$, we derive
\begin{equation}
\label{equal}
(f+\varphi)\circ\tau(x) = \varphi(x) +c\quad\text{for all }x\in\T\,.
\end{equation}

Since $\varphi$ is continuous on $\T$, equation (\ref{equal}) implies
that $(f+\varphi)\circ\tau$ is also continuous on $\T$.  In
particular, $(f+\varphi)\circ\tau$ is continuous at each discontinuity
$x$ of $\tau$.  This discontinuity is of type $(i,j)$ for some $1\le
i,j\le p$, so there exists $\varepsilon>0$ such that $\tau(z)=T_i(z)$
for $z\in(x-\varepsilon,x)$ and $\tau(z)=T_j(z)$ for
$z\in(x,x+\varepsilon)$.  So $T^{-1}(x)\cap F=\{y,y'\}$ where
$y'=T_i(x)$ is the right endpoint of a petal in $F$, and $y=T_j(x)$ is
the left endpoint of a petal in $F$.  Therefore
\begin{align*}
  (f+\varphi)(y)&=(f+\varphi)(T_j(x))=\lim_{z\nearrow x}
  (f+\varphi)(T_j(z))\cr &=\lim_{z\searrow x} (f+\varphi)(T_i(z))
  =(f+\varphi)(T_i(x))=(f+\varphi)(y')\,.
\end{align*}

Now $f+\varphi$ is Lipschitz, hence absolutely continuous, so 
by the
fundamental theorem of calculus,
\begin{align*}
  \int_y^{y'} (f+\varphi)' &= (f+\varphi)(y') -(f+\varphi)(y)=0\cr &=
  (f+\varphi)(y)-(f+\varphi)(y') = \int_{y'}^{y} (f+\varphi)' \,.
\end{align*}
In particular,
\begin{equation}
\label{ftc}
\int_{I_x} (f+\varphi)' = 0\,.
\end{equation}

Now iteration of (\ref{equal}) gives
$
\varphi =-mc + \sum_{n=1}^m f\circ\tau^n + \varphi\circ\tau^m\,,
$
and differentiation yields
$
\varphi'= \sum_{n=1}^m (f\circ\tau^n)' + (\varphi\circ\tau^m)'$
in $L^\infty$.
But $(\varphi\circ\tau^m)'=\varphi'\circ\tau^m \cdot (\tau^m)'\to0$
in $L^\infty$, because $\|(\tau^m)'\|_{L^\infty}\to0$ as $m\to\infty$
by (\ref{iteratederivbound}), so
\begin{equation}
\label{derivformula}
\varphi'= \sum_{n=1}^\infty (f\circ\tau^n)'\quad\text{in }L^\infty\,.
\end{equation}
Substituting (\ref{derivformula}) into (\ref{ftc}) gives
$
\int_{I_x} \sum_{n=0}^\infty (f\circ \tau^n)' =0
$,
as required.

(a) $\Leftarrow$ (b): By Lemma \ref{characteristicequality}, the
characteristic function $\chi(F)$ can be expressed as a linear
combination of the characteristic functions $\chi(I_x)$, for
$x\in\D$. So the integral of a function over $F$ is a linear
combination of its integrals over the intervals $I_x$, for $x\in \D$.
In particular, condition (b) implies that
$$
\int_F \sum_{n=0}^\infty (f\circ \tau^n)'=0\,,
$$
which is equivalent, by change of variable, to
\begin{equation}
\label{derivativeinterpretation}
\int \sum_{n=1}^\infty (f\circ \tau^n)'=0\,.
\end{equation}
The function $x\mapsto\sum_{n=1}^\infty (f\circ \tau^n)'(x)$, defined
Lebesgue almost everywhere, is $L^\infty$: each summand $x\mapsto
(f\circ \tau^n)'(x)$ is in $L^\infty$ because $f$ 
is Lipschitz and $\tau$ 
is piecewise
Lipschitz, and the sum converges in $L^\infty$ since $(f\circ
\tau^n)'(x)=f'(\tau^n x)(\tau^n)'(x)$, and
$\|(\tau^n)'\|_{L^\infty}\le K^{-n}$ by (\ref{iteratederivbound}).
So by (\ref{derivativeinterpretation}) the function
$x\mapsto\sum_{n=1}^\infty (f\circ \tau^n)'(x)$ lies in $L^\infty_0$,
and by Lemma \ref{rademachertheorem} it is the derivative of
some Lipschitz
continuous function $\varphi:\T\to\R$.  
Now
\begin{equation}
\label{aeonF}
(\varphi\circ T - (f+\varphi))'=0\quad\text{Lebesgue almost everywhere
on }F\,,
\end{equation}
since $\tau\circ T$ is Lebesgue almost everywhere equal to the identity
function on $F$.
But $\varphi\circ T - (f+\varphi)$ is continuous,
so by (\ref{aeonF}) its restriction to $F$ is a constant function. 

(b) $\Leftrightarrow$ (c): If $x$ is a discontinuity of $\tau$,
then
\begin{align*}
  \int_{I_x} \sum_{n=0}^\infty (f\circ \tau^n)' &=\sum_{n=0}^\infty
  \int_{I_x} (f\circ \tau^n)'=\sum_{n=0}^\infty \int_{\tau^n I_x}
  f'\cr &=\sum_{n=0}^\infty \int \chi(\tau^n I_x) f'=\int
  \sum_{n=0}^\infty \chi(\tau^n I_x) f'= \int e_x f'\,,
\end{align*}
and the equivalence of (b) and (c) follows.
\end{proof}

\begin{remark}\label{longremark}
\item[\, (a)]
  A further condition equivalent to those of Theorem
  \ref{equivalent} is that for all $x\in\D$,
\begin{equation}
\label{furtherequivalent}
\int_y^{y'} \sum_{n=0}^\infty (f\circ \tau^n)' =0
=\int_{y'}^y \sum_{n=0}^\infty (f\circ \tau^n)' \,.
\end{equation}
Clearly (\ref{furtherequivalent}) implies condition (b) of Theorem
\ref{equivalent}, while the fact that (\ref{furtherequivalent}) is
implied by condition (a) of Theorem \ref{equivalent} was
essentially established during the proof that (a) $\Rightarrow$ (b).

If, for $x\in\D$, we define $d_x:\T\to\R$ by
$$
d_x=\sum_{n=0}^\infty \chi(\tau^n J_x)\,,
$$
where
\begin{equation*}
J_x=
\begin{cases}
  [y'(x),y(x)] & \text{if }y(x)<_{y_1} y'(x)\,,\cr [y(x),y'(x)] & \text{if
  }y'(x)<_{y_1} y(x)\,,
\end{cases}
\end{equation*}
then (\ref{furtherequivalent}) becomes
$$
\int_{I_x} \sum_{n=0}^\infty (f\circ \tau^n)' =0 =\int_{J_x}
\sum_{n=0}^\infty (f\circ \tau^n)' \quad\text{for all }x\in\D\,,
$$
from which it easily follows that the condition
\begin{equation}
 \int d_x f'=0\quad \text{for all }x\in\D
\end{equation}
is also equivalent to those of Theorem \ref{equivalent}.
\item[\, (b)]
  The proof of Theorem \ref{equivalent} implies that if the
  Lipschitz function $f$ can be Lipschitz flattened on the flower $F$,
  then there is a \emph{unique} Lipschitz coboundary $g$ such that
  $(f+g)|_F$ is constant.  This is because the corresponding 
pre-image  selector $\tau$ and 
all of its iterates $\tau^n$ are uniquely
  defined on a set of full Lebesgue measure, and if we write
  $g=\varphi-\varphi\circ T$ then $\varphi$ is uniquely defined, up to
  an additive constant, by $\varphi'= \sum_{n=1}^\infty
  (f\circ\tau^n)'$ (see (\ref{derivformula})).
\end{remark}

\begin{theorem}\label{independentfunctionals}
Let $T:\T\to\T$ be an 
orientation-preserving Lipschitz expanding map.
If $F\subset\T$ is a $p$-flower then 
$\lip_F$ is a codimension-$p$ subspace of $\lip$.
\end{theorem}
\begin{proof}
Let $x_1,\ldots,x_p$ denote the discontinuities
of a pre-image selector $\tau$ for $F$. 
For $1\le j\le p$, the function $e_{x_j}$
will be denoted simply by
$e_j$. Define the linear functional $L_{F,j}:\lip\to\R$
by
$L_{F,j}(f)=\int e_j f'$.
By Theorem \ref{equivalent},
$$
\lip_F=\{f\in\lip:L_{F,j}(f)=0\text{ for }1\le j\le p\}\,,
$$
so $\lip_F$ has codimension $p$ in $\lip$ if and only if
the $p$ functionals
$\{L_{F,j}\}_{j=1}^p$ are linearly independent.
This linear independence 
is equivalent to the fact that 
$$
\int \left(\sum_{j=1}^p 
\alpha_j e_j\right) f'=0\ \text{for all }f\in\lip
\quad\Rightarrow\quad (\alpha_1,\ldots,\alpha_p)=(0,\ldots,0)\,,
$$
and since 
$f'\in L^\infty_0$ for all $f\in\lip$, this is equivalent
to the fact that
$$
\sum_{j=1}^p \alpha_j e_j\ \text{is a constant function}
\quad\Rightarrow\quad (\alpha_1,\ldots,\alpha_p)=(0,\ldots,0)\,.
$$

So suppose that
$\sum_{j=1}^p \alpha_j e_j$
is a constant function.
The strategy for showing that $\alpha_j=0$
for $1\le j\le p$
will be to consider the values of the functions
$e_j$ on the $p$ connected components of the complement of $F$,
and on the image 
under $\tau$ of a particular one of these components.

For $1\le j\le p$, define the intervals $I_j=I_{x_j}$ 
as in Notation
\ref{orderingdefn}, with reference to the counterclockwise
ordering $<_{y_1}$
(i.e.~the smallest
element is $y_1$, 
the unique pre-image of $x_1$ which is a left endpoint
  of some petal in $F$).
Denote the $p$ petals of $F$ by $P_1,\ldots, P_p$, ordered 
counterclockwise and such that $P_1$ is the petal whose left
endpoint is $y_1$.
Let $Q_1,\ldots,Q_p$ denote the $p$ connected components
of $\T\setminus F$, ordered counterclockwise and such that
$Q_1$ is the component whose left endpoint is the right endpoint
of $P_1$.

First of all note that $Q_p$ does not intersect any of the intervals
$I_{j}$. Moreover $Q_p$ does not intersect $F=\overline{\tau(\T)}$,
hence does not intersect any of the sets $\tau^n(I_{j})$
for $1\le j\le p$, $n\ge1$.
Therefore each function $e_j=\sum_{n\ge0} \chi(\tau^n I_{j})$ 
is identically zero on the interval
$Q_p$. So if $\sum_{j=1}^p \alpha_j e_j$
is a constant function on $\T$ then this constant must be zero.

It remains to show that if the function
$\sum_{j=1}^p \alpha_j e_j$ is identically zero
then $\alpha_j=0$ for $1\le j\le p$.
For this we will consider the function on the
other $p-1$ components $Q_1,\ldots,Q_{p-1}$ of the complement of $F$,
as well as on the set $\tau(Q_p)$.

First consider the restriction of the functions
$e_j$ to the
remaining $p-1$ 
components $Q_1,\ldots,Q_{p-1}$ of the complement of $F$.
Let 
$$
\b=\{j:x_j\in\A\}\quad,\quad\b'=\{j:x_j\in\A'\}\,,
$$
and for each $1\le i\le p$ define
$$
\b_i=\{j\in\b: Q_i\subset I_j\}\quad,\quad
\b_i'=\{j\in\b': Q_i\subset I_j\}\,.
$$
Note that for any $1\le i,j\le p$, 
$Q_i$ is either a subset of $I_j$ or is disjoint from $I_j$,
and
\begin{equation*}
e_j|_{Q_i} \equiv 
\begin{cases}
0&\text{when }Q_i\cap I_j=\emptyset\,,\cr
1&\text{when }Q_i\subset I_j\,.
\end{cases}
\end{equation*}
The restriction to $Q_i$
of the equation $\sum_{j=1}^p \alpha_j e_j=0$
therefore yields
\begin{equation}
\label{firstsystem}
\sum_{j\in \b_i\cup\b_i'}\alpha_j = 0\quad\text{for all }1\le i\le p-1\,.
\end{equation}
It will be convenient to consider the 
case $i=1$ of (\ref{firstsystem}), 
together with the system of equations obtained by
subtracting the $(i-1)$-st equation (\ref{firstsystem})
from
the $i$-th equation (\ref{firstsystem})
for $2\le i\le p-1$.
Since $Q_i$ is disjoint from $F$, Lemma \ref{characteristicequality}
implies that
\begin{equation}
\label{equalonQi}
\sum_{j\in\b} \chi(I_j) = \sum_{j\in\b'} \chi(I_j)
\quad\text{on }Q_i,\text{ for all }1\le i\le p\,.
\end{equation}

For each $2\le i\le p-1$, if we write
$P_i=[y_{k_i},y_{l_i}']$ then, in view of
(\ref{equalonQi}),
subtracting the $(i-1)$-st equation (\ref{firstsystem})
from
the $i$-th equation (\ref{firstsystem}) yields
\begin{equation}
\label{bothinsame}
\alpha_{k_i}=\alpha_{l_i}\quad\text{if }k_i,l_i\in\b\ \text{ or }\ k_i,l_i\in\b'\,,
\end{equation}
and
\begin{equation}
\label{inopposite}
\alpha_{k_i}=-\alpha_{l_i}\quad\text{if }k_i\in\b,l_i\in\b'\ \text{ or }\ k_i\in\b',l_i\in\b\,.
\end{equation}

Now consider $\{1,\ldots,p\}$ as the vertex set 
for an undirected graph $\Gamma$, where there is an edge
between $k$ and $l$ if and only if some petal of $F$ is
equal to $[y_k,y_l']$.
If we can prove that $\Gamma$ is connected then the equation
(\ref{firstsystem}) with $i=1$, together with the equations
(\ref{bothinsame}) and (\ref{inopposite})
for $2\le i\le p-1$, imply that
\begin{equation}
\label{1diml}
(\alpha_1,\ldots,\alpha_p)=t(\beta_1,\ldots,\beta_p)
\end{equation}
for some $t\in\R$, where 
\begin{equation*}
\beta_j=
\begin{cases}
1&\text{if $j\in\b$},\cr 
-1&\text{if $j\in\b'$.}
\end{cases}
\end{equation*}

To prove that $\Gamma$ is connected,
suppose for a contradiction that it is not, and let
$\c\subsetneq\{1,\ldots,p\}$ be the vertex set corresponding
to some connected component of $\Gamma$.
Let $\d=\{y_k\in\T:k\in\c\}\cup\{y_k'\in\T:k\in\c\}$, and let
$E\subset F$ denote the union of those petals of $F$ whose
endpoints lie in $\d$. Now
at least one petal of $F$ is disjoint from $E$,
and every petal in $E$ has positive
length, so since $T$ is Lipschitz,
\begin{equation}
\label{not0or1}
0<\text{Leb}(T(E))<1\,.
\end{equation}
For each $1\le j \le p$,
Lebesgue almost every point in $\T$ has the same
number of pre-images (under $T$) lying in $I_j$.
But
$$
\chi(E)=\sum_{j\in\c\cap\b}\chi(I_j) - \sum_{j\in\c\cap\b'}\chi(I_j)\,,
$$
so Lebesgue almost every point in $\T$ has the same
number of pre-images (under $T$) lying in $E$,
which contradicts (\ref{not0or1}). So $\Gamma$
is in fact connected, and therefore
(\ref{1diml}) holds.

It remains to show that in fact $t=0$ in (\ref{1diml}).
For this we shall consider the values of the functions $e_j$
on the set $\tau(Q_p)$.
More precisely, if $P_i$ is a petal of $F$ whose
interior has non-empty intersection with $\tau(Q_p)$ then
we shall consider the values of the functions $e_j$
on $int(P_i)\cap\tau(Q_p)$.
Note that each $e_j$ is identically equal to either
$0$ or $1$ on $int(P_i)\cap\tau(Q_p)$; it is identically
equal to $1$ if and only if $I_{j}$ contains $P_i$.

Let us write $P_i=[y_k,y_l']$.
If $j\in\{1,\ldots,p\}\setminus\{l\}$ then
$Q_i\subset I_j$ if and only if $int(P_i)\cap\tau(Q_p)\subset I_j$,
so the constant value of $e_j$ on $Q_i$ is the same as its
constant value on $int(P_i)\cap\tau(Q_p)$.
If $l\in\b$ then $Q_i\cap I_l=\emptyset$ and $P_i\subset I_l$,
hence $int(P_i)\cap\tau(Q_p)\subset I_l$,
so $e_j|_{Q_i}\equiv 0$ while
$e_j|_{int(P_i)\cap\tau(Q_p)}\equiv 1$.
If $l\in\b'$ then $Q_i\subset I_l$ and
$int(P_i)\cap I_l=\emptyset$, hence
$(int(P_i)\cap\tau(Q_p))\cap I_l=\emptyset$,
so $e_j|_{Q_i}\equiv 1$ while
$e_j|_{int(P_i)\cap\tau(Q_p)}\equiv 0$.
Therefore, subtracting the restriction of the equation
$\sum_{j=1}^p \alpha_j e_j=0$
to $int(P_i)\cap\tau(Q_p)$ 
from the restriction of the same equation
to $Q_i$ yields
$$
\alpha_l=0\,,
$$
and from (\ref{1diml}) we deduce
that $(\alpha_1,\ldots,\alpha_p)=(0,\ldots,0)$, as required.
\end{proof}

\section{Flattening on 1-flowers}
\label{1flowersection}

Let $T:\T\to\T$ 
be an orientation-preserving Lipschitz expanding map.
Let $F$ be a $1$-flower for $T$, and let $e_F$ denote the
corresponding escape time function, defined (cf.~Remark
\ref{escaperemark}) by
$$e_F= \sum_{n\ge0} \chi(\tau_F^n F)\,.$$

To prove 
Theorem \ref{flattentheorem} below it will be useful to know
that the $L^1$ function $e_F$ varies continuously with the $1$-flower
$F$.  The set of all such $1$-flowers forms a one-parameter family
$(F_\gamma)_{\gamma\in\T}$.  Since every $1$-flower $F_\gamma$ is in
particular a closed proper sub-interval of $\T$, 
with endpoints $a(\gamma)$ and $b(\gamma)$, say,
we may write
$F_\gamma=[a(\gamma),b(\gamma)]$, where both $\gamma\mapsto a(\gamma)$
and $\gamma\mapsto b(\gamma)$ are degree-one homeomorphisms of $\T$.

\begin{prop}
\label{escapecontinuous}
The map $\T\to L^1$, defined by $\gamma\mapsto e_{F_\gamma}$, is
continuous
\end{prop}
\begin{proof}
  Given any $\eta>0$, we shall show that there exists $\xi>0$ such
  that if $d(\gamma,\delta)<\xi$ then $\int 
|e_{F_\gamma} - e_{F_\delta}|  <\eta$.
  
  Let $|G|$ denote the Lebesgue measure of a measurable subset
  $G\subset\T$.  If $K>1$ is the exanding constant as in
  (\ref{expdef}) then
  $$|\tau_\gamma^nF_\gamma|\le K^{-n}|F_\gamma|< K^{-n}$$
  and
  $$
  |\tau_\delta^nF_\delta|\le K^{-n}|F_\delta|< K^{-n}$$
  by
  (\ref{iteratederivbound}), so for all $N\in\N$,
  \begin{align*}
    \int | \sum_{n>N} \chi(\tau_\gamma^nF_\gamma) - \sum_{n>N}
    \chi(\tau_\delta^nF_\delta)| &\le \sum_{n>N}
    \int\chi(\tau_\gamma^nF_\gamma) + \chi(\tau_\delta^nF_\delta) \cr
    &< 2 \sum_{n>N} K^{-n} = \frac{2K^{-(N+1)}}{1-K^{-1}}\,.
  \end{align*}
  In particular, we may choose $N$ sufficiently large so that
  $$
  \int | \sum_{n>N} \chi(\tau_\gamma^nF_\gamma) - \sum_{n>N}
  \chi(\tau_\delta^nF_\delta)| < \eta/2\,.
  $$
  
  It remains to show that we can find $\xi>0$ such that if
  $d(\gamma,\delta)<\xi$ then
\begin{equation*}
\int | \sum_{n=0}^N \chi(\tau_\gamma^nF_\gamma)
- \sum_{n=0}^N \chi(\tau_\delta^nF_\delta)|
< \eta/2\,,
\end{equation*}
or in other words,
\begin{equation}
\label{firstN}
\sum_{n=0}^N 
|\tau_\gamma^n F_\gamma\, \bigtriangleup\, \tau_\delta^n
F_\delta| < \eta/2\,,
\end{equation}
where $\bigtriangleup$ denotes symmetric difference.

The two ingredients for proving (\ref{firstN}) are that, if $\gamma$
and $\delta$ are close then firstly $F_\gamma\bigtriangleup F_\delta$
is small, and secondly the functions $\tau_\gamma$ and $\tau_\delta$
agree except on a small set. We now make this precise.

The continuity of $\gamma\mapsto a(\gamma)$ and $\gamma\mapsto
b(\gamma)$ means that for any $\varepsilon_0>0$ we may choose $\xi>0$
such that if $d(\gamma,\delta)<\xi$ then
$d(a(\gamma),a(\delta))<\varepsilon_0$ and
$d(b(\gamma),b(\delta))<\varepsilon_0$.  In particular,
\begin{equation}
\label{n=0symmetric}
|F_\gamma\bigtriangleup F_\delta|< 2\varepsilon_0\,.
\end{equation}

Let $d(\gamma,\delta)$ be small enough so that $F_\gamma$ and
$F_\delta$ intersect.  The maps $\tau_\gamma$ and $\tau_\delta$ are
identical except on the interval between their respective points of
discontinuity $T(a(\gamma))=T(b(\gamma))$ and
$T(a(\delta))=T(b(\delta))$.  If we denote this interval by
$A=A(\gamma,\delta)$, then
\begin{equation}
\label{disagreesize}
|A|=d(T(a(\gamma)),T(a(\delta))) \le C\, d(a(\gamma),a(\delta)) <
C\, \varepsilon_0\,,
\end{equation}
where $C>1$ is the Lipschitz constant for $T$
(cf.~(\ref{lipschitzT})).

Now $\tau_\gamma^n F_\gamma\, \bigtriangleup\, \tau_\delta^n F_\delta$
is contained in
\begin{align*}
  \tau_\gamma( \tau_\gamma^{n-1}F_\gamma &\cap
  \tau_\delta^{n-1}F_\delta \cap A)\ \cup\ \tau_\delta(
  \tau_\gamma^{n-1}F_\gamma \cap \tau_\delta^{n-1}F_\delta \cap A)\cr
  &\cup\ \tau_\gamma( \tau_\gamma^{n-1}F_\gamma \setminus
  \tau_\delta^{n-1}F_\delta)\ \cup\ \tau_\delta(
  \tau_\delta^{n-1}F_\delta \setminus \tau_\gamma^{n-1}F_\gamma)\,,
\end{align*}
which is itself a subset of
$$
\tau_\gamma(A)\ \cup\ \tau_\delta(A)\ \cup\ \tau_\gamma(
\tau_\gamma^{n-1}F_\gamma\, \bigtriangleup\,
\tau_\delta^{n-1}F_\delta)\ \cup\ \tau_\delta(
\tau_\gamma^{n-1}F_\gamma\, \bigtriangleup\,
\tau_\delta^{n-1}F_\delta)\,.
$$
Combining with (\ref{tauderivativebound}) and (\ref{disagreesize}),
it follows that for all $n\ge1$ we have the recurrence relation
\begin{equation}
\label{recurrencerelation}
|\tau_\gamma^n F_\gamma\, \bigtriangleup\, \tau_\delta^n F_\delta| \le
2K^{-1}\left( C\varepsilon_0 + |\tau_\gamma^{n-1} F_\gamma\, \bigtriangleup\,
  \tau_\delta^{n-1} F_\delta|\right)\,.
\end{equation}
In particular, (\ref{n=0symmetric}) and (\ref{recurrencerelation})
mean that
$$
\sum_{n=0}^N |\tau_\gamma^n F_\gamma\, \bigtriangleup\,
\tau_\delta^n F_\delta| \le B\, \varepsilon_0\,,
$$
for a constant $B=B(N,C,K)>0$ which is independent of $\gamma$ and
$\delta$.  Choosing $\varepsilon_0=\eta/(2B)$ establishes
(\ref{firstN}), as required.
\end{proof}

\begin{remark}
  For the map $T(x)=2x \pmod 1$, Bousch \cite{bousch} has given a
  quantitative bound on the modulus of continuity of the map
  $\gamma\mapsto e_{F_\gamma}$.
\end{remark}

\begin{cor}
\label{morecontinuity}
The map $\gamma\mapsto \int f' e_{F_\gamma}$ is continuous.
\end{cor}
\begin{proof}
  For any fixed $g\in L^\infty(\T)$, the 
linear functional $L^1\to\R$
  defined by $h\mapsto\int gh$ is clearly continuous, with norm
  $\|g\|_{L^\infty}$.  Now $\gamma\mapsto e_{F_\gamma}$ 
is continuous, by
  Proposition \ref{escapecontinuous}, therefore so is $\gamma\mapsto
  \int g\, e_{F_\gamma}$, and the result follows by choosing $g=f'$.
\end{proof}

We shall need 
the following well known result regarding the invariant
measures carried by $1$-flowers.

\begin{lemma}
\label{sturmianfacts}
Let $T:\T\to\T$ be an orientation-preserving expanding map.  
\item[\, (a)]
Every
$1$-flower carries 
a unique $T$-invariant probability measure; any
such measure will be called \emph{Sturmian}.
\item[\, (b)]
The support of a Sturmian
measure is equal 
to $\cap_{n\ge0} \overline{\tau^n(F)}$, where $\tau$
is either of the two pre-image selectors associated to $F$.
\item[\, (c)]
To each Sturmian measure $S$ there is an associated closed interval
$\Gamma_S$, such that the $1$-flowers carrying $S$ are precisely
$\{F_\gamma\}_{\gamma\in\Gamma_S}$, where
$F_\gamma=[a(\gamma),b(\gamma)]$. This closed interval $\Gamma_S$ is
reduced to a point if and only if $S$ is not a periodic orbit.
\end{lemma}
\begin{proof}
$T$ is topologically conjugate to $T_k(x)=kx \pmod 1$
(see \cite[p.~73, Thm.~2.4.6]{katokhasselblatt}),
where $k$ is the degree of $T$.
The conjugacy sends $1$-flowers of $T$ to $1$-flowers of $T_k$
(i.e.~closed intervals of length $1/k$), and conjugates the
corresponding pre-image selectors.
The result for $T$ therefore follows from the result for $T_k$,
and this can be proved by a straightforward adaptation of
the approach of either \cite{bouschmairesse}  or \cite{bullettsentenac}
 for the case $k=2$.
\end{proof}

\begin{remark}
\item[\, (a)]
Every Sturmian measures is ergodic, and the restriction of
$T$ to its support is combinatorially equivalent to
a rotation (see \cite{bouschmairesse,bullettsentenac}). 
More generally, the dynamics on the maximal closed
invariant subset of any flower is combinatorially equivalent
to an interval exchange transformation (see \cite{bremont1}).

\item[\, (b)]
The terminology
  Sturmian goes back to Morse \& Hedlund \cite{morsehedlund}, who
  considered certain symbol sequences on a two letter alphabet. These
  sequences correspond, under the natural symbolic coding of $T$, to
  orbits of generic points for our Sturmian measures (see
  e.g.~\cite{bousch,bouschmairesse,bullettsentenac,jenkinson2}). Some
  authors (see e.g.~\cite{lothaire, pytheasfogg}) prefer the term
  \emph{balanced} rather than Sturmian, reserving the term Sturmian
  for the \emph{non-periodic} case.

\end{remark}

By analogy with Definition \ref{flatteneddefn} we introduce the
following notions:

\begin{defn}\label{normalformmax}
  Let $T:\T\to\T$ be an expanding map.  Recall that the set $\m$ of
  $T$-invariant Borel
probability measures is compact
  for the weak$^*$ topology.  Let
  $$
  \alpha(f)=\max_{\mu\in\m} \int f\, d\mu
  $$
  denote the maximum ergodic average of the continuous function
  $f:\T\to\R$.
  
  The function $f$ is said to be in \emph{normal form} if
  $f\le\alpha(f)$.
  
  Let $F$ be a flower for $T$.  We say that $f$ can be
  \emph{continuously maximized} on $F$ if there exists a weak
  coboundary $g:\T\to\R$ such that the set of global maxima of the
  function $f+g$ is precisely $F$.  We say that $f$ can be
  \emph{Lipschitz maximized} on $F$ if there exists a Lipschitz
  coboundary $g:\T\to\R$ such that the set of global maxima of the
  function $f+g$ is precisely $F$ (or, equivalently, there exists a
  Lipschitz function $\varphi:\T\to\R$ such that the set of global
  maxima of $f+\varphi-\varphi\circ T$ is precisely $F$).
\end{defn}

\begin{remark}
\item[\, (a)]
Clearly if $f$ is continuously (respectively Lipschitz) maximized by
$F$ then it is 
continuously (respectively Lipschitz) flattened by $F$.
Moreover, 
if $g$ is the corresponding weak coboundary then $f+g$ is in
normal form (since $F$ carries a $T$-invariant probability measure, 
so $\alpha(f)=\max f$); 
therefore the $f$-maximizing measures are
precisely those carried by $F$.
\item[\, (b)]
If $F$ is a $1$-flower then the notion of a function being maximized
on $F$ corresponds to the \emph{Sturmian condition} defined by Bousch
\cite{bousch}.  
\item[\, (c)]
Br\'emont \cite{bremont2} considers Lipschitz
functions maximized on flowers, showing that such functions can be
Lipschitz approximated by functions with periodic maximizing measures.
\end{remark}

The following well known result guarantees that if $f$ is Lipschitz
then we can add a Lipschitz coboundary to it so that the resulting
function is in normal form.

\begin{lemma}\label{normalform}
  Let $T:\T\to\T$ be an expanding map, and let $f:\T\to\R$ be
  Lipschitz. There exists a Lipschitz function $\varphi$ such that
  $f+\varphi-\varphi\circ T\le \alpha(f)$.
\end{lemma}
\begin{proof}
  This result seems to date 
back to an unpublished manuscript of Conze
  \& Guivarc'h \cite{conzeguivarch}.  Published proofs can be found in
  \cite[Thm.~1]{bouschwalters} and \cite[Thm.~4.7]{lectures}, while
  the proofs of \cite[Lem.~A]{bousch} and
  \cite[Thm.~9]{contreraslopesthieullen}, where the stated hypotheses
  are slightly stronger than ours, are also easily adapted.
\end{proof}

\begin{theorem}
\label{flattentheorem}
Let $T:\T\to\T$ be an orientation-preserving 
Lipschitz expanding map, and let $f:\T\to\R$ be
Lipschitz.  If a Sturmian measure $S\in\m$ is $f$-maximizing, then
there exists a 1-flower $F$
carrying $S$ such that $f$ can be Lipschitz flattened on $F$.
\end{theorem}

\begin{proof}
  In view of Lemma \ref{normalform}, it will suffice to show that if
  $f\le\alpha(f)$ then $f$ can be Lipschitz flattened on $F$.  So
  assume that the Lipschitz function $f$ satisfies $f\le\alpha(f)$,
  and the Sturmian measure $S$ is $f$-maximizing. Note that such an
  $f$ is identically equal to $\alpha(f)$ on $\text{supp}(S)$.

  By Theorem \ref{equivalent}
  it suffices to show that there exists a 1-flower $F$ containing
  $\text{supp}(S)$ such that
\begin{equation}\label{vanishes}
\int f' e_F = 0\,,
\end{equation}
where
$$
e_F = \sum_{n\ge0} \chi(\tau_F^n F)
$$
is the escape time function for $F$.

Let $F_{\gamma^-}=[a(\gamma^-),b(\gamma^-)]$ denote the 1-flower whose
right endpoint is the rightmost point of $\text{supp}(S)$, and
$F_{\gamma^+}=[a(\gamma^+),b(\gamma^+)]$ the 1-flower whose left
endpoint is the leftmost point of $\text{supp}(S)$.  Of course if $S$
is a non-periodic Sturmian measure then $F_{\gamma^-}=F_{\gamma^+}$ is
the unique 1-flower containing $\text{supp}(S)$, while if $S$ is
periodic then $F_{\gamma^-}\neq F_{\gamma^+}$.

We claim that
\begin{equation}
\label{posneg}
\int f' e_{F_{\gamma^-}}  \ge 0
\quad\text{and}\quad
\int f' e_{F_{\gamma^+}}  \le 0\,.
\end{equation}

Now $\gamma\mapsto \int f' e_{F_\gamma}$ is continuous by Corollary
\ref{morecontinuity}, so once (\ref{posneg}) is established, the
intermediate value theorem will imply the existence of
$\gamma_0\in[\gamma^-,\gamma^+]$ such that (\ref{vanishes}) holds for
$F=F_{\gamma_0}=[a(\gamma_0),b(\gamma_0)]$.

So to prove the theorem it remains to prove (\ref{posneg}).  In fact
we shall only prove that $\int f' e_{F_{\gamma^-}} \ge 0$, as the
proof of the other inequality is analogous.  To simplify notation we
shall write $G=F_{\gamma^-}$.  We have
\begin{equation}
\label{integral}
\int f' e_{G} =\int f' \sum_{n\ge0} \chi(\tau_G^nG) 
=\sum_{n\ge0} \int_{\tau_G^n G} f'\,.
\end{equation}

Now $G$ is an interval, and $\tau_G$ has a single jump discontinuity,
so each $\tau_G^nG$ is a union of $s_n\le n+1$ intervals $I^{(n)}_j$,
which we write as
$$
\tau_G^nG= \cup_{j=1}^{s_n} I^{(n)}_j\,.
$$
Let $c_j^{(n)}$ and $d_j^{(n)}$ denote, respectively, the left and
right endpoints of $I_j^{(n)}$.  The right endpoint of $G$ is a point
in $\text{supp}(S)$, therefore every right endpoint $d_j^{(n)}$ is
also a point in $\text{supp}(S)$.  It follows that
\begin{equation}
\label{dominate}
f(d_j^{(n)})= \alpha(f)=\max f \ge f(c_j^{(n)})\,.
\end{equation}
Now $f$ is Lipschitz, hence absolutely continuous, so applying the
fundamental theorem of calculus, and then (\ref{dominate}), gives
\begin{equation*}
\int_{c_j^{(n)}}^{d_j^{(n)}} f' = 
f(d_j^{(n)}) - f(c_j^{(n)})\ge0\,.
\end{equation*}

Therefore for all $n\ge0$,
\begin{equation*}
  \int_{\tau_G^n G} f'
= \sum_{i=1}^{s_n} \int_{c_j^{(n)}}^{d_j^{(n)}} f'
\ge0\,,
\end{equation*}
so from (\ref{integral}) we deduce that
$$\int f' e_{G} =\sum_{n\ge0} \int_{\tau_G^n G} f' \ge0\,,
$$
as required.
\end{proof}

\begin{remark}
\label{zerofinding}
As mentioned in Section \ref{introsection}, for certain functions $f$
the maximizing measure is known to be Sturmian, but a priori it is
not known \emph{which} of the Sturmian measures is maximizing.
A consequence of
Theorem \ref{flattentheorem} is that in order to show that a 
particular Sturmian measure $S$ is $f$-maximizing, it suffices
to locate a $1$-flower $F$ which carries $S$ and on which $f$
can be Lipschitz
flattened. If $\gamma\mapsto F_\gamma$ is a parametrisation
of the $1$-flowers for $T$ then, by Theorem \ref{equivalent}, we
must find $\gamma$ such that 
\begin{equation}
\label{eqnnumerical}
\int e_{F_\gamma}f'=0\,.
\end{equation}
The equation (\ref{eqnnumerical}) can be solved 
numerically\footnote{An approximate solution is typically 
sufficient, in view of the fact that periodic Sturmian measures
are carried by a parameter interval of $1$-flowers.}
by approximating the escape time functions $e_{F_\gamma}
=\sum_{n=0}^\infty \chi(\tau_{F_\gamma}^n(F_\gamma))$
by finite truncations 
$e_{F_\gamma,N}:=\sum_{n=0}^N \chi(\tau_{F_\gamma}^n(F_\gamma))$.
The distance $\|e_{F_\gamma}-e_{F_\gamma,N}\|_{L^1}$
decreases exponentially with $N$, so solutions
to the equation $\int e_{F_\gamma,N}f'=0$, which can be computed
using a root-finding algorithm such as Newton's
method, converge to solutions of (\ref{eqnnumerical})
at an exponential rate. 
\end{remark}

\begin{question}\label{unresolved}
  Does some analogue of Theorem \ref{flattentheorem} hold for more
  general flowers?  
For example suppose $\mu$ is the unique invariant
  measure carried by some flower $F$, and is the unique maximizing
  measure for a Lipschitz function $f$. Is it then the case that $f$
  can be Lipschitz flattened on some flower (not necessarily $F$)
  which carries $\mu$ (and which necessarily carries no other
  invariant measure)?  The proof of any such result would seem to
  require a higher dimensional analogue of the intermediate value
  theorem.
  
  This question also raises the issue of whether or not every
  invariant measure carried by \emph{some} flower is in fact the
  unique invariant measure carried by some (other) flower.
\end{question}

\begin{example}\label{notmax}
  If a Lipschitz function $f$ has a Sturmian maximizing measure $S$,
  Theorem \ref{flattentheorem} guarantees that $f$ can be Lipschitz
  flattened on some $1$-flower which carries $S$.  One might expect
  that in fact such a 1-flower Lipschitz \emph{maximizes} $f$, in the
  sense of Definition \ref{normalformmax} (e.g.~this is exactly what
  Bousch \cite{bousch} proves in the case where $T(x)=2x \pmod 1$ and
  $f$ is a trigonometric polynomial of degree one).  However in
  general this is not the case: for example if $f$ has a Sturmian
  maximizing measure, but this is not the unique $f$-maximizing
  measure, then clearly no $1$-flower can maximize $f$.  In fact even
  when the Sturmian measure is the \emph{unique} maximizing measure,
  it is not the case that there exists a $1$-flower which maximizes
  $f$, as the following example shows.
  
  Consider the expanding map $T(x)=2x \pmod 1$, whose $1$-flowers are
  precisely the set of all closed semi-circles (i.e.~intervals of
  length $1/2$).  There are infinitely many $\gamma\in(0,1/6)$ such
  that the $1$-flower $F_\gamma=[\gamma,\gamma+1/2]$ contains a
  non-periodic Sturmian measure $S$ (see e.g.~\cite{bullettsentenac}).
  Fix one such $\gamma$, and let $\tau$ denote the pre-image selector
  for $F_\gamma$ defined by
\begin{equation*}
\tau(x)=
\begin{cases}
  \frac{x+1}{2}&\text{ if }x\in[0,2\gamma)\cr \frac{x}{2}&\text{ if
  }x\in[2\gamma,1)\,.
\end{cases}
\end{equation*}
Since $\text{supp}(S)$ is the intersection of the decreasing sequence
$\overline{\tau^n(F_\gamma)}$, the Sturmian measure $S$ is in
particular carried by
$$\overline{\tau(F_\gamma)} =[\gamma,\gamma/2+1/4] \ \cup\ 
[\gamma/2+1/2,\gamma+1/2]\,.$$

Let $f:\T\to\R$ be the continuous piecewise linear function whose
maximum value is $0$ and whose derivative is given by
\begin{equation*}
f'\equiv
\begin{cases}
  0&\text{ on }(\gamma,\gamma/2+1/4)\cr -2/\gamma & \text{ on
  }(\gamma/2 +1/4, \gamma+1/4)\cr 2/(1/2-\gamma) &\text{ on
  }(\gamma+1/4,\gamma/2+1/2)\cr 0 &\text{ on
  }(\gamma/2+1/2,\gamma+1/2)\cr -1&\text{ on
  }(\gamma+1/2,\gamma+3/4)\cr 1&\text{ on }(\gamma+3/4,\gamma)\,.
\end{cases}
\end{equation*} 

Note that $f$ is in normal form: its set of global maxima $f^{-1}(0)
=[\gamma,\frac{\gamma}{2}+\frac{1}{4}] \ \cup\ 
[\frac{\gamma}{2}+\frac{1}{2},\gamma+\frac{1}{2}]$ carries the
Sturmian measure $S$, and this is the unique maximizing measure.  Now
$F_\gamma$ is the only $1$-flower which contains $S$, so if $f$ is
Lipschitz maximized by a $1$-flower then it must be Lipschitz
maximized by $F_\gamma$.
There is a unique Lipschitz coboundary $g$ such that
$(f+g)|_{F_\gamma}$ is a constant (cf.~Remark \ref{longremark}), so to
show that $f$ is \emph{not} Lipschitz maximized by $F_\gamma$ it
suffices to show that $F_\gamma$ is not the set of maxima of the
function $f+g$.

Let $\varphi:\T\to\R$ be a continuous piecewise linear function
(uniquely defined up to an additive constant) whose derivative is
given by
\begin{equation*}
\varphi'\equiv
\begin{cases}
  0&\text{ on }(\gamma,\gamma+1/2)\cr -1/\gamma&\text{ on
  }(\gamma+1/2,2\gamma+1/2)\cr 1/(1/2-\gamma)&\text{ on
  }(2\gamma+1/2,\gamma)\,,
\end{cases}
\end{equation*}
and define the Lipschitz coboundary $g$ by
$$g=\varphi-\varphi\circ T\,.$$
It is readily verified that
$$(f+g)|_{F_\gamma}\equiv 0\,.$$
However,
\begin{align*}
  (f&+g)(\gamma+3/4) =(f+g)(\gamma+3/4)-(f+g)(\gamma+1/2) \cr&=
  \int_{\gamma+1/2}^{\gamma+3/4} (f+g)'\cr &=
  \int_{\gamma+1/2}^{2\gamma+1/2} (f+g)'
  +\int_{2\gamma+1/2}^{\gamma/2+3/4} (f+g)'
  +\int_{\gamma/2+3/4}^{\gamma+3/4} (f+g)'\cr
  &=-\gamma\left(1+\frac{1}{\gamma}\right)+\left(\frac{1}{4}-
    \frac{3\gamma}{2}\right)\left(\frac{1}{1/2-\gamma}-1\right)
  +\frac{\gamma}{2}\left(\frac{1}{1/2-\gamma}+\frac{2}{\gamma}-1\right)\cr
  &= \frac{1}{4} - \frac{\gamma}{1-2\gamma}\ ,
\end{align*}
which is strictly positive because $\gamma\in(0,1/6)$.  So $0$ is not
the maximum value of $f+g$, and therefore $F_\gamma$ is not the set of
maxima of $f+g$.  Therefore $F_\gamma$ does not Lipschitz maximize
$f$, as required.
\end{example}

\WARMprocessEPS{notmax_b}{eps}{bb}
\begin{figure}[htpb]
$$\begin{xy}
\xyMarkedImport{}
        \xyMarkedPos{1}*!L\txt\labeltextstyle{(a) $f$}
        \xyMarkedPos{2}*!L\txt\labeltextstyle{(b) $g = \varphi-\varphi\circ T$}
        \xyMarkedPos{3}*!L\txt\labeltextstyle{(c) $f+g$}
        \xyMarkedPos{4}*\txt\labeltextstyle{\tiny{$\gamma$}}
        \xyMarkedPos{5}*\txt\labeltextstyle{\tiny{$\frac{\gamma}{2} +$\textonequarter}}
        \xyMarkedPos{6}*\txt\labeltextstyle{\tiny{$\gamma + $\textonequarter}}
        \xyMarkedPos{7}*\txt\labeltextstyle{\tiny{$\frac{\gamma}{2} + $\textonehalf}}
        \xyMarkedPos{8}*\txt\labeltextstyle{\tiny{$\gamma +$\textonehalf}}
        \xyMarkedPos{23}*\txt\labeltextstyle{\tiny{$\gamma +$\textthreequarters}}
        \xyMarkedPos{11}*!R\txt\labeltextstyle{}
        \xyMarkedPos{9}*\txt\labeltextstyle{\tiny{$-1$}}
        \xyMarkedPos{10}*\txt\labeltextstyle{\tiny{$0$}}
        \xyMarkedPos{12}*\txt\labeltextstyle{\tiny{$-1$}}
        \xyMarkedPos{16}*\txt\labeltextstyle{\tiny{$0$}}
        \xyMarkedPos{13}*\txt\labeltextstyle{\tiny{$-1$}}
        \xyMarkedPos{22}*\txt\labeltextstyle{\tiny{$0$}}
\end{xy}$$
\caption[Non-maximizing function]{Flattening of
the function $f$ from Example \ref{notmax}: 
$f$ itself is shown in (a), 
the unique Lipschitz coboundary $g = \varphi-\varphi\circ T$
flattening $f$ on the one-flower
$F= [\lambda, \lambda+1/2]$ is shown in (b), and
the flat function $f+g$ is
shown in (c).  Although flat on $F$, the function $f+g$
is not in normal form as its maximum (at $\gamma + 3/4$) is larger than $0$.
The function $f$ is in normal form;
the unique Sturmian measure carried by $F$ is also
carried by the set of maxima of $f$, so is the unique
$f$-maximizing measure. 
Thus $f$ may be flattened on $F$
(i.e.~satisfies the pre-Sturmian condition),
but the corresponding flat function is not in normal form,
so $f$ does not satisfy the Sturmian condition,
even though its maximizing measure is Sturmian.}
\label{fig:notmax}

\end{figure}
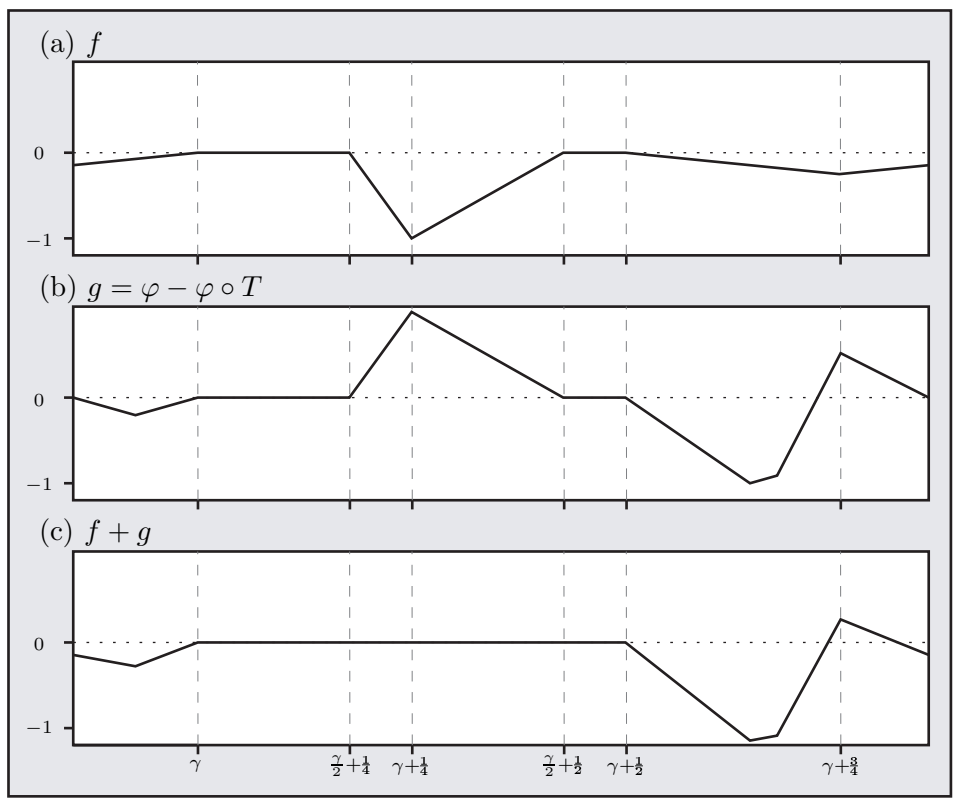

\end{document}